\documentclass{amsart}
\usepackage{amsmath, amssymb, amscd}

\newtheorem{thm}{Theorem}
\newenvironment{example}{\vspace{.1in} \noindent {\bf Example}:}{\vspace{.3in}}
\newenvironment{Definition}{\vspace{.1in} \noindent {\bf Definition}:}{\vspace{.3in}}

\newcommand{\cL}{{\mathcal L}}
\newcommand{\fM}{{\mathfrak M}}
\newcommand{\fX}{{\mathfrak X}}

\begin{document}

\title{Uniformity in the Mordell-Lang conjecture}
\author{Thomas Scanlon}
\thanks{Partially supported by NSF Grant DMS-0071890}
\email{scanlon@math.berkeley.edu}
\address{University of California, Berkeley \\
Department of Mathematics \\
Evans Hall \\
Berkeley, CA 94720-3840 \\
USA}
\date{24 May 2001}
\maketitle

\section*{Introduction}
The Mordell-Lang conjecture and its variants assert 
subvarieties of algebraic groups can meet certain 
subgroups only in a finite set of translates of subgroups.
In most cases the number of such translates is unknown
and it is not even known whether this number may be bounded
by a function of the geometric data.  In this note we 
show that some uniformity follows immediately from the finiteness
result.   

The main technical result behind this note is a theorem of 
Pillay on the stability of the theory of an algebraically 
closed field with a predicate for a group of Lang type~\cite{Pillay}.  My 
interest in these questions was renewed through my reading
of the paper~\cite{Remond} in which a version of uniformity for the 
Mordell-Lang conjecture for abelian varieties is proved.  
As the reader will see the main result of this note is an 
immediate corollary of Pillay's theorem, but somehow it was not 
noticed earlier.

\section*{Main Theorem}
We begin by stating more precisely what we mean by 
an analogue of the Mordell-Lang Conjecture.

\begin{Definition}
We say that the subgroup $\Gamma$ of an algebraic group $G$
over an algebraically closed field $K$ is of \emph{Lang type} if
\begin{itemize}
\item for any subvariety $X$ of a Cartesian 
	power $G^n$ of $G$ the set $X(K) \cap \Gamma^n$ is a finite 
	union of cosets of subgroups.
\end{itemize}
We say that $\Gamma$ is \ \emph{uniformly of Lang type} if in 
addition

\begin{itemize}
\item for any algebraic family $\{ X_b \}_{b \in B}$ of subvarieties
of $G^n$ there is a sequence $H_1, \ldots, H_m$ of algebraic 
subgroups of $G$ such that for any $b \in B$ there are a 
set $I \subseteq \{1, \ldots, m\}$ and points $\gamma_1, 
\ldots, \gamma_m \in \Gamma^n$ such that the Zariski closure of
$X_b(K) \cap \Gamma^n$ is $\bigcup_{i \in I} \gamma_i + H_i$.
\end{itemize} 
\end{Definition}

There are a number of known groups of Lang type.

\begin{example}
\begin{itemize}
\item If $G$ is a semiabelian variety defined over a field $K$ of
	characteristic zero and $\Gamma \leq G(K)$ is 
	a subgroup of finite dimension ($\dim_{\mathbb Q} \Gamma 
	\otimes {\mathbb Q} < \infty$), then $\Gamma$ is 
	of Lang type~\cite{Faltings,Vojta,McQuillan}.
\item If $A$ is an abelian variety over a field $K$ of 
	characteristic $p$ admitting no non-trivial 
	homomorphisms of algebraic groups to abelian varieties
	defined over a finite field and $\Gamma \leq A(K)$
	is a subgroup of finite ${\mathbb Z}_p$-rank 
	(${\rm rk}_{{\mathbb Z}_p} \Gamma \otimes {\mathbb Z}_p
	< \infty$), then $\Gamma$ is of Lang type~\cite{Hrushovski}.

\item If $K$ is a field of positive characteristic and 
	$\Gamma \leq K$ is torsion module of a Drinfeld 
	module of generic characteristic, then $\Gamma$ is of Lang type~\cite{Scanlon}.
\end{itemize} 
\end{example}

Recall that a group $G$ considered as an $\cL$-structure for 
some first-order language $\cL$ is said to be \emph{modular} if for 
every natural number $n$ every $\cL$-definable subset of $G^n$ is 
a finite Boolean combination of cosets definable subgroups. 

Recall that if $\fM$ is an $\cL$-structure for some first-order
language $\cL$ and $X \subseteq M^n$ is a nonempty subset of some 
Cartesian power of $M$, the the full-induced structure on $X$ is 
the structure $\fX$ having universe $X$ and basic relations 
$D \cap X^m$ for each $\cL_M$-definable $D \subseteq M^{nm}$.

We could restate the definition of $\Gamma$ being of Lang type 
as \emph{the full induced structure
on $\Gamma$ from $(K,+,\cdot)$ is modular}.

\begin{thm}
If $\Gamma$ is of Lang type,  then it is uniformly of Lang type.
\end{thm}
\begin{proof}
By Proposition 2.6 of~\cite{Pillay}, the structure $(K, +, \cdot,
\Gamma)$ is stable and the formula $x \in \Gamma$ is one-based. 
Suppose that $\Gamma$ is not uniformly of Lang type 
witnessed by some family $\{ X_b \}_{b \in B}$ of 
subvarieties of $G^n$.  Then for any 
finite sequence $B_1, \ldots, B_m$ of subgroups of $\Gamma^n$
there is some $b \in B$ for which the formula 
$(\forall x_1, \ldots, x_m \in G^n) \bigwedge_{I \subseteq \{ 1, \ldots,
m \} } (\exists y \in \Gamma^n) [(y \in X \land \bigwedge_{i \in I} y \notin 
x_i + B_i) \vee (y \notin X \land \bigvee_{i \in I} y \in x_i + B_i)]$ holds.
By compactness, we can find some elementary extensions $({^*}K, +, \cdot, {^*}\Gamma)$
and a point $b \in B({^*}K)$ satisfying all of these formulas for
all choices of finite sequences of subgroups of $\Gamma^n$.  As $\Gamma$ 
is one-based as a definable group in $(K, +, \cdot, \Gamma)$, ${^*}\Gamma$ is 
also one-based. Thus, every (parametrically) definable 
(in $({^*}K, +, \cdot, {^*}\Gamma)$) subset of ${^*}\Gamma^n$ is 
a finite Boolean combination of cosets of $(K, +, \cdot, \Gamma)$-definable
subgroups.  Thus, $X_b({^*}K) \cap {^*}\Gamma$, being a definable subset of
${^*}\Gamma^n$, is a finite Boolean combination of cosets of
$(K, +, \cdot, \Gamma)$-definable
subgroups of ${^*}\Gamma^m$.  As $X_b$ is Zariski closed, this combination 
is necessarily a finite union.  Thus, there are subgroups $B_1, \ldots, B_m$
of $\Gamma$ and points $a_1, \ldots, a_m \in G({^*}K)$ such that 
$({^*}K, +, \cdot, {^*}\Gamma) \models 
(\forall y \in {^*}\Gamma^n) (y \in X_b \leftrightarrow \bigvee_{i = 1}^m 
y \in a_i + B_i)$.  This contradicts the choice of $b$.

Therefore, $\Gamma$ is uniformly of Lang type.
\end{proof}


\begin{thebibliography}{99}

\bibitem{Faltings} {\sc G. Faltings}, Diophantine approximation 
on abelian varieties,  \emph{Ann. of Math. (2)} {\bf 133}
(1991), no. 3, 549--576.

\bibitem{Hrushovski} {\sc E. Hrushovski}, The Mordell-Lang
 conjecture for function fields, \emph{J. Amer. Math. Soc.} {\bf 9}
(1996), no. 3, 667--690. 

\bibitem{McQuillan} {\sc M. McQuillan}, Division points 
on semi-abelian varieties, \emph{Invent. Math.}
 {\bf 120} (1995), no. 1, 143--159. 



\bibitem{Pillay} {\sc A. Pillay}, The model-theoretic content
of Lang's conjecture, {\bf Model theory and algebraic geometry},
{\bf LNM 1696}, Springer, Berlin, 1998, 101--106.

\bibitem{Remond} {\sc G. R\'{e}mond}, D\'{e}compte dans une
conjecture de Lang, \emph{Invent. Math.} {\bf 142} (2000), no. 3, 
513--545.

\bibitem{Scanlon} {\sc T. Scanlon}, The Manin-Mumford conjecture for
 Drinfeld modules, preprint, 1999.

\bibitem{Vojta} {\sc P. Vojta}, Siegel's theorem in the compact 
case, \emph{Ann. of Math. (2)} {\bf 133} (1991), no. 3, 509--548. 

\end{thebibliography}
\end{document}